\documentclass[11pt]{amsart}
\usepackage{amscd,amssymb,amsmath}
\usepackage[all]{xy}

\newtheorem{theorem}{Theorem}[section]
\newtheorem{lemma}[theorem]{Lemma}
\theoremstyle{definition}

\newtheorem{proposition}[theorem]{Proposition}
\newtheorem{corollary}[theorem]{Corollary}
\newtheorem{conjecture}[theorem]{Conjecture}
\newtheorem{question}[theorem]{Question}

\newtheorem{fact}{Fact}  
\theoremstyle{remark}
\newtheorem{remark}{Remark}[section]
\newcommand{\til}[1]{\widetilde{#1}}
\newcommand{\be}{\begin{enumerate}}
\newcommand{\ee}{\end{enumerate}}
\newcommand{\bq}{\begin{question}}
\newcommand{\eq}{\end{question}}
\newcommand{\bcj}{\begin{conjecture}}
\newcommand{\ecj}{\end{conjecture}}
\newcommand{\bc}{\begin{corollary}}
\newcommand{\ec}{\end{corollary}}
\newcommand{\bl}{\begin{lemma}}
\newcommand{\el}{\end{lemma}}
\newcommand{\btl}{\begin{technicalLemma}}
\newcommand{\etl}{\end{technicalLemma}}
\newcommand{\bt}{\begin{theorem}}
\newcommand{\et}{\end{theorem}}
\newcommand{\bp}{\begin{proposition}}
\newcommand{\ep}{\end{proposition}}
\newcommand{\bft}{\begin{fact}}
\newcommand{\eft}{\end{fact}}
\newcommand{\brk}{\begin{remark}}
\newcommand{\erk}{\end{remark}}
\newcommand{\bd}{\begin{Dn}}
\newcommand{\ed}{\end{Dn}}
\newcommand{\Mat}[4]{\left( \begin{array}{rr}
                            #1 & #2 \\
                            #3 & #4
                      \end{array} \right)}

\def\t{\tau}

\def\t{\tau}

\def\Id{\operatorname{Id}}

\newcommand{\ov}[1]{\overline{#1}}
\newcommand{\coker}{{\rm Coker\ }}
\numberwithin{equation}{section}
\newcommand{\bea} {\begin{eqnarray*}}
\newcommand{\beq} {\begin{equation}}
\newcommand{\bey} {\begin{eqnarray}}
\newcommand{\eea} {\end{eqnarray*}}
\newcommand{\eeq} {\end{equation}}
\newcommand{\eey} {\end{eqnarray}}

\newcommand{\F}{$F$ } 
\newcommand{\Z}{ \mathbf Z }

\newcommand{\R}{\mathbf R }

\begin{document}
\title[Twisted conjugacy classes in R. Thompson's group \F]{Twisted conjugacy classes in R. Thompson's group F}
\author{Collin Bleak}
\address{ Department of Mathematics\\
Cornell University\\ 
Ithaca, NY 14853-4201, USA\\}
\email{bleak@math.cornell.edu}
\author{Alexander Fel'shtyn}
\address{ Instytut Matematyki, Uniwersytet Szczecinski,
ul. Wielkopolska 15, 70-451 Szczecin, Poland
and Boise State University, 1910
University Drive, Boise, Idaho, 83725-155, USA }
\email{felshtyn@diamond.boisestate.edu, felshtyn@mpim-bonn.mpg.de}
\author {Daciberg L. Gon\c{c}alves}
\address{Dept. de Matem\'atica - IME - USP, Caixa Postal 66.281 - CEP 05311-970,
S\~ao Paulo - SP, Brasil}
\email{dlgoncal@ime.usp.br}

\begin{abstract}
In this short article, we prove that any automorphism of the R. Thompson's group $F$ has   infinitely many twisted conjugacy classes.  The result follows from the work of Brin  \cite{BrinCh}, together with a standard facts  on R. Thompson's group $F$, and elementary properties of the Reidemeister numbers.
\end{abstract}

\date{\today}
\keywords{Reidemeister number, twisted conjugacy classes, R. Thompson's   group \F,
 $R_\infty$ property}
\subjclass{20E45;37C25; 55M20}
\maketitle

\bibliography{ref}
\bibliographystyle{amsplain}

\section{Introduction}

Let $\phi:G\to G$ be an automorphism of a group $G$. A class of equivalence
$x\sim gx\phi(g^{-1})$ is called the \emph{Reidemeister class }
or $\phi$-\emph{conjugacy class} or \emph{twisted conjugacy class of} $\phi$. The
number $R(\phi)$ of Reidemeister classes is called the \emph{Reidemeister number}
of $\phi$.
The interest in twisted conjugacy relations has its origins, in particular,
in the Nielsen-Reidemeister fixed point theory (see, e.g. \cite{j, FelshB}),
in Selberg theory (see, eg. \cite{Shokra,Arthur}), and  Algebraic Geometry
(see, e.g. \cite{Groth}).

A current important  problem of the field concerns obtaining 
a twisted analogue of the celebrated Burnside-Frobenius theorem
\cite{FelHill,FelshB,FelTro,FelTroVer,ncrmkwb,polyc,FelTroObzo},
that is, to show the
coincidence of
the Reidemeister number of $\phi$ and the number of fixed points of the
induced homeomorphism of an appropriate dual object.
One step in this process is to describe the class of groups $G$, such that $R(\phi)=\infty$ for
any automorphism $\phi:G\to G$.

The work of discovering which groups belong to the mentioned class of groups was begun by Fel'shtyn and Hill in \cite{FelHill}.
Later, it was shown by various authors that the following groups belong to this class:
 (1) non-elementary Gromov hyperbolic groups \cite{FelPOMI,ll}, (2)
Baumslag-Solitar groups $BS(m,n) = \langle a,b | ba^mb^{-1} = a^n \rangle$
except for $BS(1,1)$ \cite{FelGon},
 (3) generalized Baumslag-Solitar groups, that is, finitely generated groups
which act on a tree with all edge and vertex stabilizers infinite cyclic
\cite{LevittBaums}, (4)
lamplighter groups $\Z_n \wr \Z$ iff $2|n$ or $3|n$ \cite{gowon1}, (5)
the solvable generalization $\Gamma$ of $BS(1,n)$ given by the short exact sequence
$1 \rightarrow \Z[\frac{1}{n}] \rightarrow \Gamma \rightarrow \Z^k \rightarrow 1$
 as well as any group quasi-isometric to $\Gamma$ \cite{TabWong},
groups which are quasi-isometric to $BS(1,n)$ \cite{TabWong2} (while this property is
 not a quasi-isometry invariant),  (6) saturated weakly branch groups (including the Grigorchuk group and the Gupta-Sidki group) \cite{FelLeoTro}.

The paper \cite{TabWong} suggests a terminology for this property, which we would like to follow.  Namely, a group $G$ has \emph{property } $R_\infty$ if all of its automorphisms $\phi$
have $R(\phi)=\infty$.

For the immediate consequences of the $R_\infty$ property in topological fixed point
theory see, e.g., \cite{TabWong2}.

In the present note we prove that  R. Thompson's group \F   has the $R_\infty$ property.
We do not know if this is the case for injective homomorphisms.
The R. Thompson  group \F is a finitely-presented group which  has exponential growth, does not contain a free nonabelian subgroup; \F
is not a residually finite group
and is not an elementary amenable group \cite{CFP}. It is still unknown whether or not \F is amenable.

Richard Thompson introduced the groups $F \leq T \leq V $ \cite{Thom} in connection
with his studies in logic. Thompson's groups have since appeared in a variety of mathematical topics: the world problem, infinite simple groups, homotopy and shape theory, group cohomology, dynamical systems and analysis.   We are interested not only in Thompson's group $F$, but also in the groups $T$ and $V$, and the generalizations $F_n\leq T_n\leq V_n$ (Note that Higman introduces the groups $V_n$ and solves the conjugacy problem for these groups in \cite{HigmanFPSG}, while Brown carries this generalization out for $T_n$ and $V_n$ in \cite{BrownFinite}).  Our work here only discusses the group $F$.

The results of the present paper demonstrate that the further study
of Reidemeister theory for R. Thompson's  group \F  has to go
along the lines specific for the infinite case. On the other
hand, this result shrinks the class of groups for which the
twisted Burnside-Frobenius conjecture \cite{FelHill,FelTro,FelTroVer,ncrmkwb,polyc,FelTroObzo}
has yet to be verified.
 We would like to complete the introduction with the following question

\begin{question}
Do the generalised Thompson's groups $F_n$ have the $R_\infty$ property?
\end{question} 

\bigskip

{\bf Acknowledgments}: The first author would like to thank Professor  Brin for helpful conversations and  discussion of his work in \cite{BrinCh}, upon which this note strongly relies.  The second   author would like to thank
 V. Guba, R. Grigorchuk,   M. Kapovich and M. Sapir,   for stimulating discussions and comments. 

\section{ R. Thompson's group \F : definitions and standard facts}
We will be working with R. Thompson's group $F$.  Following the definition in \cite{BrinCh}, $F$ consists of a restricted class of homeomorphisms of the real line under the operation of composition.  A homeomorphism $\alpha: \R \to \R$ is in $F$ if and only if $\alpha$

\be
\item is piecewise-linear (admitting finitely many breaks in slope),
\item is orientation-preserving,
\item has all slopes of affine portions of its graph in the set $\{2^k|k 
\in \Z  \}$,
\item has all breaks in slope occuring over the dyadic rationals $\Z[1/2]$,
\item maps $\Z[1/2]$ into $\Z[1/2]$, and 
\item has first and last affine components equal to pure translations by (potentially distinct) integers $\alpha_l$ and $\alpha_r$.
\ee
Note that the fifth condition is actually equivalent to saying that $\alpha$ maps $\Z[1/2]$ bijectively onto $\Z[1/2]$ in an order preserving fashion, given the other axioms.

Every normal subgroup of $F$ contains the commutator subgroup $F' = [F,F]$ of $F$ (see Theorem 4.3 in the standard introductory survey \cite{CFP} on $F$). Further, in our realization of $F$, the commutator subgroup $F'$ of $F$ consists of those elements of $F$ which are the identity function near $\pm \infty$.  In particular, there is the standard onto homomorphism $Ab:F\to\Z\times\Z$, with kernel $F'$, which is defined by the rule $Ab(f) = (f_l,f_r)$ (this is Theorem 4.1 of \cite{CFP}). Here, $f_l$ is the translational part of $f$ near $-\infty$, and $f_r$ is the translational part of $f$ near $\infty$, as indicated by the similar notation in the definition above for $F$.  In fact, given any $k\in F$, let us fix the notation $k_l$ and $k_r$ as the translational parts of $k$ (near $-\infty$ and $\infty$ respectively) for the rest of this short paper.

Finally, we will need a deep result of Matthew Brin from \cite{BrinCh}.  To state this result, let us first define $Rev:F\to F$ to be the automorphism produced by conjugating an element of $F$ by the real homeomorphism $x\mapsto -x$.  We also need to define the set of \emph{eventually $T$-like} piecewise linear self-homeomorphisms (admitting an infinite discrete set of ``breaks" in the domain, where the first derivative is not defined) of $\R$ as follows.

A homeomorphism $\alpha: \R \to \R$ is \emph{eventually $T$-like} if and only if $\alpha$

\be
\item is piecewise-linear (admitting a discrete (possibly infinite) collection of breaks in slope),
\item is orientation-preserving,
\item has all slopes of affine portions of its graph in the set $\{2^k|k 
\in \Z  \}$,
\item has all breaks in slope occuring over the dyadic rationals $\Z[1/2]$,
\item maps $\Z[1/2]$ into $\Z[1/2]$, and 
\item has minimal domain value $R_{\alpha}\geq 0$ and maximal domain value $L_{\alpha}\leq 0$ so that $\alpha(x+1) = \alpha(x) +1$ for all $x\in\R\backslash(L_{\alpha}-1,R_{\alpha})$.
\ee

The final condition says that to the far left and the far right, any such $\alpha$ (which could be in $F$ in every other way) projects to circle maps, as in those locations it satisfies the periodicity equation $\alpha(x+1) = \alpha(x) +1$.  If $\alpha$ satisfied the periodicity equation across all of $\R$, then it would represent a lift of an element of R. Thompson's group $T$.  Note that elements of $F$ are eventually $T$-like.  For any eventually $T$-like element $\alpha$, we will use the notation $L_{\alpha}$ and $R_{\alpha}$ as in the above definition.

In \cite{BrinCh}, Brin shows that the inner automorphism of $Homeo(\R)$ defined by conjugating elements of $Homeo(\R)$ by any specific eventually $T$-like element restricts to an automorphism of $F$.  We will call such an automorphism of $F$ an \emph{ eventually $T$-like conjugation} of $F$.  The following is a restatement of a portion of Theorem 1 in \cite{BrinCh} using our language.

\bt(Brin)
\label{AutoClassification}
The automorphism group of $F$ is generated by $Rev$ and the eventually $T$-like conjugations of $F$.
\et

We will measure the impact of an automorphism in $\Z\times\Z$, the abelianization of $F$.  In order to do this, we will use the following lemma.

\begin{lemma}
Let $f\in F$, and let $g:\R\to\R$ be eventually $T$-like.  If $k = f^g$, then $f_l = k_l$ and $f_r = k_r$.  That is, eventually $T$-like conjugation does not change the translational parts of $f$ near plus or minus infinity.  
\end{lemma}
\begin{proof}
Let $f\in F$, and $g:\R\to\R$ be an eventually $T$-like homeomorphism of $\R$.  Let $R_g$ and $L_g$ be the locations so that for all $x\in \R\backslash (L_g-1,R_g)$ we have that  $g(x+1) = g(x) +1$.

Define variation functions $Var_f$, $Var_g:\R\to\R$ by the rules $Var_f(x) = f(x) - x$, $Var_g(x) = g(x)-x$.  Define $V_R = Var(R_g)$.  For $x>R_g$, we have that $|Var_g(x) -V_R|<1$.  That is, on a global scale, homeomorphisms satisfying the periodicity equation are close to pure translations.  In particular, the graph of $g$ is not only periodic near $\infty$, it is also always within a distance of $1$ to the pure translation by the constant $V_R$.  Note also that for $x>R_f$, $Var_f(x) = f_r$.

Let $M = |R_g| + |R_f| + |V_R| + |f_r| + 1$, and suppose $x\in[M,\infty)$.  We see that the following equations are true.
 
\[g(f(g^{-1}(x)) = g(g^{-1}(x) + f_r) = g(g^{-1}(x))+ f_r = x+f_r= f(x)\]
\vskip.15 in

(In the third to last equality, if $f_r$ is negative, we are actually using the fact that $g(y) = g(y+1)-1$, where $y > R_g$.) 

The argument near minus infinity is similar.
\end{proof}
Let  $H_1 = H_1^{Gp}$ be the first integral homology functor
from groups to abelian groups.

The above work produces this corollary.
 \begin{corollary}
 {\it If $\phi:F \to F$ is an eventually $T$-like automorphism of $F$ then the induced automorphism on the abelianization $H_1(F) \cong \Z\times\Z$ of $F$ is the identity.}
\end {corollary}

In particular, if $H_1(Rev)$ is the induced  homomorphism on the abelianization $H_1(F)$ of $F$ we have the following.

 \begin{corollary}
 {\it The image of $Aut(F)$ in  $Aut(H_1(F)=\Z\times \Z)$ is a cyclic group $\Z_2$ having $H_1(Rev)$  as a generator}.
\end{corollary}

\section{Simple facts about Reidemeister classes and the Main theorem}
\subsection{Reidemeister classes and inner automorphisms}

Let us denote by $\t_g:G\to G$ the automorphism $\t_g(x)=g x\,g^{-1}$
for $g\in G$. Its restriction on a normal subgroup we will denote by $\t_g$
as well. We will need the following statements.

\begin{lemma}\label{lem:redklassed}
$\{g\}_\phi k=\{g\,k\}_{\t_{k^{-1}}\circ\phi}$.
\end{lemma}

\begin{proof}
Let $g'=f\,g\,\phi(f^{-1})$ be $\phi$-conjugate to $g$. Then
$$
g'\,k=f\,g\,\phi(f^{-1})\,k=f\,g\,k\,k^{-1}\,\phi(f^{-1})\,k
=f\,(g\,k)\,(\t_{k^{-1}}\circ\phi)(f^{-1}).
$$
Conversely, if $g'$ is $\t_{k^{-1}}\circ\phi$-conjugate to $g$, then
$$
g'\,k^{-1}=
f\,g\,(\t_{k^{-1}}\circ\phi)(f^{-1})k^{-1}=
f\,g\,k^{-1}\,\phi(f^{-1}).
$$
Hence a shift maps $\phi$-conjugacy classes onto classes related to
another automorphism.
\end{proof}

\begin{corollary}\label{lem:innerreidem}
$R(\phi)=R(\t_g \circ \phi)$. In particular $R(\Id)=R(\t_g)$ is the number of the usual conjugacy classes in $G$.
\end{corollary}

Consider a group extension respecting homomorphism $\phi$:
\begin{equation}\label{eq:extens}
 \xymatrix{
0\ar[r]&
H \ar[r]^i \ar[d]_{\phi'}&  G\ar[r]^p \ar[d]^{\phi} & G/H \ar[d]^{\ov{\phi}}
\ar[r]&0\\
0\ar[r]&H\ar[r]^i & G\ar[r]^p &G/H\ar[r]& 0,}
\end{equation}
where $H$ is a normal subgroup of $G$.
First, notice that the Reidemeister classes of $\phi$ in $G$
are mapped epimorphically onto classes of $\ov\phi$ in $G/H$. Indeed,
\begin{equation}\label{eq:epiofclassforexs}
p(\til g) p(g) \ov\phi(p(\til g^{-1}))= p (\til g g \phi(\til g^{-1}).
\end{equation}
Suppose that the Reidemeister number  $R(\ov\phi)$ is infinite, the previous remark then implies that the Reidemeister number
$R(\phi)$ is infinite.

See \cite{go:nil1} for a generalization of this elementary fact  to  homomorphisms of short exact sequences.
An endomorphism $\phi:G\rightarrow G$ is said to be eventually
commutative if there exists a natural number $n$ such that the
subgroup $\phi^n(G)$ is commutative.

We are now ready to compare the Reidemeister number of an
endomorphism $\phi$ with the Reidemeister number
of $H_1(\phi):H_1(G)\rightarrow H_1(G)$.

\begin{theorem}[\cite{j}]
The composition $\eta\circ\theta$,
$$
 G  \stackrel{\theta}{\longrightarrow} H_1(G)
       \stackrel{\eta}{\longrightarrow}
       \coker\left[H_1(G) \stackrel{1-H_1(\phi)}{\longrightarrow} H_1(G)\right] ,
$$
where $\theta$ is abelianization and $\eta$ is the natural projection,
sends every $\phi$-conjugacy class to a single element.
Moreover, any group homomorphism $\zeta:G \rightarrow \Gamma$ which
sends every $\phi$-conjugacy class to a single element,
factors through $\eta\circ\theta$.
If $\phi:G\rightarrow G$ is eventually commutative (for example, if the  group $G$  is abelian)
then
$$
R(\phi)
 = R(H_1(\phi))
 = \#\coker (1-H_1(\phi)).
$$
\end{theorem}
The first part of this theorem  is trivial. If $\alpha^\prime=\gamma \alpha \phi (\gamma^{-1})$ , then
$$
\theta(\alpha^\prime )=\theta(\gamma) + \theta(\alpha) + \theta(\phi (\gamma^{-1}))=
$$
$$
=\theta(\gamma) + \theta(\alpha) - H_1(\phi)(\theta(\gamma))=\theta(\alpha) +(1- (H_1(\phi))\theta(\gamma),
$$
hence  $\eta\circ \theta(\alpha)= \eta\circ \theta(\alpha^\prime)$.

This Theorem shows the importance of the group $ \coker(1 -(H_1(\phi))$
\begin{corollary}(\cite{j}, p.33)

Let  $G$ be  a finitely generated  free Abelian group.
 
If $ det(I - \phi)=0$
 then $ R(\phi)=\#\coker (1-\phi)=\infty $.
\end{corollary}

We can now state the main result.

\begin{theorem}
For any automorphism $\phi$ of R.Thompson's  group \F the Reidemeister number $R(\phi)$
is infinite.

\end{theorem}

\noindent
\begin{proof}
  We will  calculate the matrix  $M$ of the automorphism $H_1(Rev)$.
 Let  $g$, $k\in F$ with $Rev(g) = k$. Then $ k_r=-g_l$ and $ k_l=-g_r$. So,  the matrix  of an  automorphism $H_1(Rev)$ is the matrix
$M=\Mat 0{-1}{-1}0$. Hence  $det( I - M )=0.$
 From corollary 3.4 it follows  that  the number of Reidemeister classes of 
$H_1(Rev)$ is infinite. The same is true  for all powers $M^k$, since for k even $M^k=id$ and $M^k=M$ for k odd.
Hence, by Corollary 2.3 and  2.4, the Reidemeister number  $R(H_1(\phi))$ is also  infinite. From (3.1)- (3.2) it follows that 
the Reidemeister number $R(\phi)$ is infinite as well.
\end{proof}

In a more concrete fashion, and for the curious reader, it is not difficult to generate an infinitude of twisted conjugacy classes in $\Z\times\Z$ for $H_1(Rev)$.  Consider the set $\Gamma$ of pairs $\left\{(0,a)|a\in\Z\right\}$,  one may check directly that no two elements of $\Gamma$ are twisted-conjugate equivalent.

\bibliographystyle{amsplain}

\end{document}